\documentclass[12pt]{amsart}
\usepackage{amssymb}
\usepackage{layout}
\numberwithin{equation}{section}
\newtheorem{theorem}{Theorem}[section]

\theoremstyle{definition}
\newtheorem{definition}[theorem]{Definition}

\theoremstyle{remark}

\newcommand{\lambdaB}{\mbox{\boldmath$\lambda$}}
\newcommand{\PhiB}{\mbox{\boldmath$\Phi$}}
\newcommand{\sigmaB}{\mbox{\boldmath$\sigma$}}
\newcommand{\tauB}{\boldsymbol{\tau}}
\newcommand{\reals}{\mathbb R}

\newcommand{\charv}[1]{Q^{#1}}
\newcommand{\comp}{\mkern-1mu\raise 1pt 
\hbox{$\,\scriptstyle\circ\,$}}
\newcommand{\crosss}[1]{K^{(#1)}}
\newcommand{\cpx}[1]{\mbox{\sffamily{k}}^{(#1)}}
\newcommand{\cpxo}[1]{\mbox{\sffamily{k}}_o^{(#1)}}
\newcommand{\cpxos}[1]{\mbox{\scriptsize{\sffamily{k}$^{(#1)}_o$}}}

\newcommand{\D}{\mathcal D}
\newcommand{\DEder}[1]{\mbox{D}_{#1}}
\newcommand{\DD}[1]{\mathbb D_{#1}}
\newcommand{\DE}[1]{\mbox{D}_{#1}}
\newcommand{\Di}{\mathcal D^{(\infty)}}

\newcommand{\Dx}[1]{\mathcal D^{(#1)}}
\newcommand{\dom}{\mbox{dom}\,}

\newcommand{\dop}[1]{\mathfrak d_{#1}}
\newcommand{\DT}[1]{\mathbf{D}_{#1}}
\newcommand{\Dxf}[2]{\mathcal D^{(#1)}{}_{\mkern-6mu\;|#2}}

\newcommand{\Ex}[1]{\mathcal E^{(#1)}}
\newcommand{\Ei}{\mathcal E^{(\infty)}}
\newcommand{\emult}[1]{\mathcal M_{#1}{}}
\newcommand\fdim{q}
\newcommand{\flow}[2]{\PhiB^{#1}_{#2}}
\newcommand{\g}{\mathfrak g}
\newcommand{\gx}[1]{J^{#1}\mathfrak g}
\newcommand{\gxf}[2]{J^{#1}_{#2}\mkern2mu\mathfrak g}%
\newcommand{\gxfo}[2]{J^{#1}_{o,#2}\mkern2mu\mathfrak g}%

\newcommand{\gind}{k}
\newcommand{\G}{\mathcal G}
\newcommand{\Gx}[1]{\mathcal G^{(#1)}}

\newcommand{\gpx}[1]{\mbox{\sffamily{g}}^{(#1)}}\newcommand{\gpxo}[1]{\mbox{\sffamily{g}}^{(#1)}_o}
\newcommand{\gpxs}[1]{\mbox{\scriptsize{\sffamily{g}$^{(#1)}$}}}

\newcommand{\hpx}[1]{\mbox{\sffamily{h}}^{(#1)}}
\newcommand{\hpxo}[1]{\mbox{\sffamily{h}}_o^{(#1)}}
\newcommand{\hpxs}[1]{\mbox{\scriptsize{\sffamily{h}$^{(#1)}$}}}
\newcommand{\hpxos}[1]{\mbox{\scriptsize{\sffamily{h}$^{(#1)}_o$}}}
\newcommand{\Hx}[1]{\mathcal H^{(#1)}}

\newcommand{\ISD}[2]{\mathcal D^{(#1)}_{\,#2}}
\newcommand{\IS}[2]{\mathcal G^{(#1)}_{\,#2}}

\newcommand{\IDE}[3]{L_{#1}{}^{#3}_{,#2}}
\newcommand{\IDEn}{L^{(\jo)}}
\newcommand{\id}{\mbox{id}}
\newcommand{\idjet}[2]{\mathbb I^{(#1)}_{#2}}
\newcommand{\jpx}[1]{\mbox{\sffamily{z}}^{(#1)}}
\newcommand{\jpxo}[1]{\text{\sffamily{z}}^{(#1)}_o}
\newcommand{\jpxs}[1]{\mbox{\scriptsize{\sffamily{z}$^{(#1)}$}}}%
\newcommand{\jpxos}[1]{\mbox{\scriptsize{\sffamily{z}}$^{(#1)}_o$}}
\newcommand{\jet}[2]{j^{#1}_{#2}}
\newcommand{\J}[1]{J^{#1}}

\newcommand{\Jxf}[2]{J^{(#1)}_{\;|#2}}
\newcommand{\Ji}{J^\infty}
\newcommand\jo{n}
\newcommand{\JTMx}[1]{J^{#1}TM}
\newcommand{\JTMxfo}[2]{J^{#1}_{o,#2}TM}
\newcommand{\JTMxf}[2]{J^{#1}_{#2}TM}

\newcommand{\lactionD}[1]{\mathcal L_{#1}}
\newcommand{\lactionE}[1]{\widehat{\mathcal L}_{#1}}
\newcommand{\lactionJ}[1]{\widetilde{\mathcal L}_{#1}}
\newcommand{\lift}[2]{\lambdaB^{(#1)}(#2)}
\newcommand{\liftmap}[1]{\lambdaB^{(#1)}}
\newcommand{\liftmapj}[2]{\lambdaB^{(#1)}_{#2}}
\newcommand{\lm}{\mbox{\L}}
\newcommand{\locdif}{\varphi}
\newcommand{\locdiff}{\psi}

\newcommand\mdim{m}
\newcommand{\mfx}[1]{\rho^{(#1)}}
\newcommand{\mpnt}{\mbox{\sffamily{z}}}

\newcommand{\mpnts}{\mbox\scriptsize{\textsf{z}}}
\newcommand{\mpntos}{\mbox\scriptsize{\textsf{z}}_o}
\newcommand{\mPnt}{\mbox{\sffamily{Z}}}%

\newcommand{\ops}{n^*}

\newcommand{\pion}{\widetilde{\pi}^n_o}
\newcommand{\projD}[2]{{\pi}^{#1}_{#2}}
\newcommand{\projE}[2]{\widehat{\pi}^{#1}_{#2}}
\newcommand{\projJ}[2]{\widetilde{\pi}^{#1}_{#2}}
\newcommand{\prc}[2]{\widehat{\phi}^{#1}_{#2}}
\newcommand{\prm}{\mbox{\textbf{pr}}}
\newcommand{\prx}[1]{\mbox{\textbf{pr}}^{(#1)}}
\newcommand{\prvx}[2]{\phi^{#1}_{#2}}

\newcommand{\pro}[1]{\mbox{pr}^{(#1)}\,}

\newcommand{\ractionD}[1]{\mathcal R_{#1}}
\newcommand{\range}{\mbox{im}\;}

\newcommand\sdim{p}
\newcommand{\smfld}{\mathcal R}
\newcommand{\smfldm}{\mathcal S}
\newcommand{\sourceDx}[1]{\sigmaB^{#1}}
\newcommand{\sourceEx}[1]{\widehat\sigmaB^{#1}}
\newcommand{\sourceHx}[1]{\widehat\sigmaB^{#1}_{\mathcal H}}
\newcommand{\sups}[2]{#1^1,\ldots,#1^{#2}}

\newcommand{\targetDx}[1]{\tauB^{#1}}
\newcommand{\targetEx}[1]{\widehat\tauB^{#1}}
\newcommand{\targetHx}[1]{\widehat\tauB^{#1}_{\mathcal H}}
\newcommand{\targetHCx}[1]{\mu^{(#1)}}
\newcommand{\targetHCinvx}[1]{\eta^{(#1)}}

\newcommand{\Uset}{\mathcal U}
\newcommand{\Usetx}[1]{\mathcal U^{(#1)}}

\newcommand{\Uhat}[2]{\widehat{U}^{#1}_{#2}}

\newcommand{\Un}{U^{(n)}}

\newcommand{\uc}[2]{u^{#1}_{#2}}
\newcommand{\uco}[2]{u^{#1}_{o,#2}}
\newcommand{\Uc}[2]{U^{#1}_{#2}}

\newcommand{\vb}{\mathbf v}

\newcommand{\vpntn}{\mbox{\sffamily{v}}^{(n)}}

\newcommand{\vpntxo}[1]{\mbox{\sffamily{v}}^{(#1)}_o}
\newcommand{\Vset}[1]{\mathcal V^{#1}}
\newcommand{\VsetT}[1]{\widetilde{\mathcal V}^{#1}}
\newcommand{\where}{\mbox{where}}
\newcommand{\wpntx}[1]{\mbox{\sffamily{w}}^{(#1)}}
\newcommand{\wpntxhat}[1]{\widehat{\mbox{\sffamily{w}}}^{(#1)}}
\newcommand{\W}[2]{W^{#1}_{#2}}
\newcommand{\Wset}{\mathcal W}

\newcommand{\vcx}[2]{\zeta^{#1}_{#2}}
\newcommand{\vxx}[1]{\xi^{#1}}
\newcommand{\vux}[1]{\phi^{#1}}
\newcommand{\zx}[1]{z^{(#1)}}
\newcommand{\zc}[1]{z^{#1}}
\newcommand{\zoc}[1]{z_o^{#1}}

\newcommand{\Zc}[2]{Z{}^{#1}_{#2}}
\newcommand{\Zoc}[2]{Z_{\lower1.1pt\hbox{$\scriptstyle o,$}}{}^{#1}_{#2}}
\newcommand{\Zn}{Z^{(n)}}

\newcommand{\zZc}[1]{(z,Z^{(#1)})}
\newcommand{\zZn}{(z,Z^{(n)})}
\newcommand{\zcfx}[2]{\zeta^{#1}_{#2}}
\newcommand{\zxc}[2]{\vxx{#1}_{#2}}
\newcommand{\zuc}[2]{\vux{#1}_{#2}}
\newcommand{\zcfxo}[2]{\zeta_{\lower1.1pt\hbox{$\scriptstyle o,$}}{}^{#1}_{#2}}
\newcommand{\zzcfn}{(z,\zeta^{(n)})}

\newcommand{\ZEn}{\mathbf Z^{(n)}}
\newcommand{\zngn}{(\jpx n,\gpx n)}
\newcommand{\znZn}{(\zx n,\Zn)}

\newcommand{\xc}[1]{x^{#1}}
\newcommand{\xu}[1]{(x,u^{(#1)})}
\newcommand{\X}{\mathcal X}
\newcommand{\XDx}[1]{\X(\Dx{#1})}
\newcommand{\XRISD}[2]{\mathcal X_R^{}(\mathcal D^{(#1)}_{#2})}
\newcommand{\Xc}[2]{X^{#1}_{#2}}
\newcommand{\XM}{\mathcal X(M)}
\newcommand{\Xn}{X^{(n)}}

\begin{document}

\title[Persistence of Freeness]{Persistence of Freeness for \\ Lie Pseudogroup Actions}

\author[Peter Olver]{Peter J. Olver$^\dagger$}
\address{School of Mathematics, University of Minnesota, Minneapolis, MN 55455, U.S.A.}
\thanks{$^\dagger$ Supported in part by NSF Grant 08--07317.}
\email{olver@math.umn.edu}
\urladdr{http://www.math.umn.edu/\~{}olver/}

\author{Juha Pohjanpelto}
\address{Department of Mathematics, Oregon State University, Corvallis, OR 97331, U.S.A.}
\email{juha@math.oregonstate.edu}
\urladdr{http://oregonstate.edu/\~{}pohjanpp/}

\textwidth=400pt

\begin{abstract}{The action of a Lie pseudogroup $\G$ on a smooth manifold $M$ induces a prolonged pseudogroup action on the jet spaces $\J\jo$ of submanifolds of $M$. 
We prove in this paper that both the local and global freeness of the action of $\G$ on $\J\jo$ persist under prolongation in the jet order $\jo$. Our results underlie the construction of complete moving frames and, indirectly, their applications in the identification and analysis of the various invariant objects for the pseudogroup action on $\J\infty$}.
\end{abstract}



\maketitle

\section{Introduction}
The results in this paper are motivated by recent developments in the study of pseudo-groups, their moving frames and invariants, and a range of applications, \cite{COP2,OP2,OP3,EDS2008}.
The classical treatments \cite{Cartan1937, Green1978, Griffiths1974, Jensen1977} of moving frames are primarily concerned with equivalence, symmetry and rigidity properties of submanifolds $S\subset G/H$ of homogeneous spaces under the natural action of $G$. Moving frames in these time-honored problems may in effect be identified as suitably normalized equivariant local sections on $S$ of the bundle $G\to G/H$, or as lifts to $G$ of maps into $G/H$ by means of such sections. A more general point of view is adopted in \cite{FelsOlver1999}, where an alternate description of a moving frame is put forth as an equivariant section of the action groupoid $M\times G\to M$ associated with the Lie group action on a manifold $M$.  This reformulation served to open a wide range of applications reaching well beyond those afforded by the classical approach to moving frames.  See \cite{Olversurvey} for a recent survey of activity in this area.

Given an infinite dimensional pseudogroup $\G$ acting on $M$, our main focus lies on its induced action on submanifolds $\smfldm\subset M$. In this framework the principal protagonists are the jet spaces $\Gx\jo$ of pseudogroup transformations and $\J\jo$ of submanifolds of $M$, endowed with the natural action of $\G$.
For pseudogroups, which are
characterized via their action on a manifold, the proper analogue of the finite dimensional action groupoid is furnished by the bundles $\Ex\jo\to\J\jo$ composed of pairs $(\jpx\jo,\gpx\jo)$ of jets $\jpx\jo\in\J\jo$ and $\gpx\jo\in\Gx\jo$ with the same base point in $M$. Moving frames can then be conceived as local sections of $\Ex\jo$ equivariant under the joint action of $\G$ on the constituent spaces, and are, likewise to the finite dimensional situation, ordinarily constructed via a normalization process based on a choice of a cross section to the pseudogroup orbits in $\J\jo$, cf.~\cite{OP2}.   

In concrete applications one frequently deals with moving frames of increasingly high order that are mutually compatible under the natural projections $\projJ{\jo+\gind}{\jo}\colon\J{\jo+\gind}\to\J\jo$. These, by the way of projective limits, collectively form a so-called complete moving frame on $\J\infty$. 
As expounded in \cite{OP2,OP3}, complete moving frames, when combined with Gr\"obner basis techniques,
can be effectively used to identify differential invariants, invariant differential forms, operators of invariant differentiation, and so on, for the prolonged action of Lie pseudogroups on $\J\infty$, and to uncover the algebraic structure of the invariants and of the invariant variational bicomplex, \cite{KOivb}. We refer to \cite{COP1, COP2,Moro2005,EDS2008,ShemMans2008} for recent applications involving the method of moving frames for infinite dimensional pseudogroups.

In the finite dimensional situation of a Lie group action, the existence of a moving frame requires that the action be locally free, \cite{FelsOlver1999}. However, as bona fide infinite dimensional groups cannot have trivial isotropy, one is lead to define (local) freeness of the action in terms of the jets of group transformations fixing a point in $\J\jo$, \cite{OP2}. The adapted definition relying on jets constrains the dimensions of the jet spaces $\Gx\jo$, and 
provide a simpler alternative to the Spencer cohomological growth conditions imposed by Kumpera, \cite{Kumpera}, in his
analysis of differential invariants.  
Our notion of freeness, when applied to finite dimensional group actions, proves to be slightly broader than the classical concept,
and, as we will elaborate in section \ref{S:MF} of the present paper, ensures the existence of local moving frames for pseudogroup actions on $\J\jo$. By contrast, extending the moving frame method and results to non-free actions remains an open problem.

Since freeness is the essential attribute in our constructions, our first order of business is to establish its persistence under prolongations. Specifically, as the main contributions of the present paper, we prove in Theorems \ref{T:LocalFreeness} and \ref{T:GlobalFreeness} that if a pseudogroup acts (locally) freely at $\jpx\jo\in\J\jo$, then it also acts (locally) freely at any $\jpx{\jo+\gind}\in\J{\jo+\gind}$, $\gind\geq0$, 
with $\projJ{\jo+\gind}{\jo}(\jpx{\jo+\gind})=\jpx\jo$.
These results, notably, are the key ingredients to the construction of complete moving frames and, indirectly, underlie the various applications requiring invariant quantities for pseudogroup actions and the analysis of their algebraic structure. The local result, Theorem \ref{T:LocalFreeness}, appeared in its original form in \cite{OP3}, with a proof resting on techniques from commutative algebra. Here we give an alternate, direct proof of the Theorem requiring only basic linear algebra. The global result of Theorem \ref{T:GlobalFreeness} is new and highlights the differences between the classical finite dimensional theory of group actions, \cite{OlverEIS}, and the infinite dimensional theory as developed in \cite{OP2}.

Our paper is organized as follows. We start in section \ref{S:CP} with an overview of the nuts and bolts of continuous pseudogroups, which is followed by an outline of prolonged pseudogroup actions on submanifold jet bundles in section \ref{S:JB}. Then, in section \ref{S:MF}, we review the method of moving frames for pseudogroup actions on submanifold jet bundles $\J\jo$. These appear in two guises --- as locally and globally equivariant sections of the bundle $\Ex\jo\to\J\jo$ associated with the prolonged action --- and we discuss conditions guaranteeing the existence of each type. Finally, in section \ref{S:PF}, we establish the main results of this paper, namely, the persistence of both local and global freeness of pseudogroup actions under prolongation in the jet order.

\section{Lie Pseudogroups}\label{S:CP}

Let $M$ be a smooth $\mdim$-dimensional manifold. 
Denote the pseudogroup of all 
local diffeomorphisms $\locdif$ of $M$ with an open domain $\dom\locdif\subset M$ by $\D = \D(M)$ and the bundle of their $\jo^{\mbox{\scriptsize th}}$ order jets $\gpx\jo=\jet\jo\mpnts\locdif$, $\mpnts\in\dom\locdif$, by $\Dx n = \Dx n(M)$. 
Write 
\begin{equation}\label{E:projD}
\projD\gind\jo\colon\Dx\gind \longrightarrow \Dx\jo,\qquad 0\leq \jo\leq\gind,
\end{equation}
for the canonical projections.
The source 
$\sourceDx{n}\colon\Dx n\to M$ 
and target maps 
$\targetDx n{}\colon\Dx n\to M$ 
are given by
\begin{equation}
\label{sourcetarget}
\sourceDx\jo({\jet\jo\mpnts\locdif}) = \mpnt,\qquad
\targetDx\jo({\jet\jo\mpnts\locdif}) = \locdif(\mpnt),
\end{equation}
respectively. Let $\Dxf\jo{\mpnts}=(\sourceDx n)^{-1}(\mpnt)$ stand for the 
source fiber and $\ISD\jo{\mpnts} = (\sourceDx n)^{-1}(\mpnt)\cap(\targetDx n)^{-1}(\mpnt)$
for the Lie group of isotropy jets at $\mpnt$, the latter being isomorphic with the \emph{prolonged general linear group}, that is, the Lie group of $\jo$-jets of local diffeomorphisms of $\reals^\mdim$ fixing the origin.  

The bundle $\Dx\jo$ is equipped with a groupoid multiplication, \cite{Mackenzie},
induced by composition of mappings,
\begin{equation}\label{E:Dmult}
\jet\jo{\locdif(\mpnts)}\locdiff\cdot\jet\jo\mpnts{\locdif}=
\jet\jo\mpnts(\locdiff\comp\locdif),\qquad\text{ $\locdif(\mpnt)\in\dom\locdiff$.}
\end{equation}
The operation \eqref{E:Dmult} also defines the actions
\begin{equation}\label{eqn:actionD}
\lactionD\locdif\gpx\jo=\jet\jo{\targetDx\jo(\gpxs\jo)}\locdif\cdot\gpx\jo,\qquad
\ractionD\locdif\gpx\jo = \gpx\jo\mkern-3mu\cdot\jet\jo\mpnts\locdif,
\end{equation}
of $\D$ on $\Dx n$ by left and right multiplication in an obvious fashion.

Given local coordinates
$\zc{}=(\sups{\zc{}{}}\mdim)$, 
$\Zc{}{}=(\sups{\Zc{}{}{}}\mdim)$
on $M$ about $\mpnt$ and $\mPnt = \locdif(\mpnt)$,
respectively, 
the induced local coordinates of 
$\gpx\jo = \jet\jo\mpnts\locdif\in\Dx\jo$ 
are given by $\zZn$, where the components 
\begin{equation}\label{E:TargetDerVars}
\Zc a{b_1b_2\cdots b_\gind}=\dfrac{\partial^{\gind}\locdif^a}{\partial \zc{b_1}\partial\zc{b_2}\cdots\partial\zc{b_\gind}}(\mpnt),   
\quad\mbox{$1\leq a\leq\mdim$,\quad $0\leq \gind\leq\jo$,}
\end{equation}
of $\Zn$ represent the partial derivatives of the
coordinate expression $\locdif^a=\Zc a{}\comp\locdif$ evaluated at the 
source point $\mpnt = \sourceDx n(\gpx n)$.
Following Cartan, 
we will use lower case letters, 
$z$, $x$, $u$, \dots for the source coordinates 
and the corresponding upper case letters 
$\Zn$, $\Xn$, $\Un$, \dots  for the derivative
target coordinates of the diffeomorphism jet $\gpx n$. 

Let $\XM$ denote the sheaf of locally defined smooth vector fields on $M$, and write $\JTMx\jo$ for the space of their $n^\text{th}$ order jets. 
Given local coordinates $\zc{}=(\sups{\zc{}{}} m)$ on  $M$, a vector field is written in component form as
\begin{equation}\label{E:VectComp} 
\vb = \sum_{a=1}^\mdim \vcx a{}(\zc{}) 
\dfrac{\partial}{\partial \zc a}\,,
\end{equation}
and the coordinates on $\JTMx n$ induced by \eqref{E:VectComp} are designated by 
\begin{equation}\label{E:TMJetCoord}
\zzcfn=(\zc a,\zcfx b{},\zcfx b{c_1},\dots,\zcfx b{c_1c_2\cdots c_n}).
\end{equation}

A vector field $\vb\in\XM$ lifts to a right-invariant vector field 
\begin{equation}\label{E:Lift}
\lift n{\vb}\in\XDx\jo
\end{equation}
defined on $(\targetDx\jo)^{-1}(\dom\vb)\subset\Dx\jo$ as the infinitesimal generator of the left action of its flow map $\flow\vb t$ on $\Dx\jo$, cf.~\cite{OP1}.
The lift $\lift\jo\vb$ is vertical, that is, tangent to the source fibers $\Dxf\jo{\mpnts}$ and has the expression
\begin{equation}\label{E:LiftMapCoord}
\lift\jo\vb = \sum_{a=1}^\mdim\sum_{\gind=0}^\jo\DD {\zc {b_1}}\DD {\zc {b_2}}\cdots\DD{\zc{b_\gind}}\vcx a{}(\Zc{}{})\dfrac{\partial}{\partial \Zc a{b_1b_2\cdots b_\gind}}
\end{equation}
in the local coordinates \eqref{E:TargetDerVars}, where
\begin{equation}\label{E:DTotDerOp}
\DD {\zc b} = \dfrac{\partial}{\partial \zc b}
+\Zc c{b}\dfrac{\partial}{\partial\Zc c{}}
+\Zc c{bc_1}\dfrac{\partial}{\partial \Zc c{c_1}}
+\Zc c{bc_1c_2}\dfrac{\partial}{\partial \Zc c{c_1c_2}}+\cdots
\end{equation}
denotes the standard coordinate total derivative operators on $\Di$. The lift map $\liftmap\jo$ is easily seen to respect the Lie brackets of vector fields. 

As is well known, the space $\JTMxfo\jo\mpnts$ of $n$-jets at $\mpnt$ of vector fields vanishing at $\mpnt$ becomes a Lie algebra when equipped with the bilinear operation induced by the usual Lie bracket of vector fields. With this operation, the lift map \eqref{E:LiftMapPoint} can be seen to restrict to an isomorphism
\begin{equation}\label{E:LiftMapPoint}
\liftmapj\jo{\mpnts}\colon \JTMxfo\jo\mpnts\>\longrightarrow\> \XRISD\jo{\mpnts},\qquad \mpnt\in M,
\end{equation}
between  $\JTMxfo\jo\mpnts$ and the Lie algebra of right-invariant vector fields on the isotropy subgroup $\displaystyle \ISD \jo{\mpnts}$.

Recall that an $n+1$ jet $\jet{n+1}\mpnts\sigma$ defines a linear map
\begin{displaymath}
\lm_{\jet{n+1}\mpnts\sigma}:T_{\mpnts} M\>\longrightarrow\> T_{\jet n\mpnts\sigma}\Dx n \qquad {\rm by} \qquad \lm_{\jet{n+1}\mpnts\sigma}\vb = (\jet{n}{}\sigma)_*\vb.
\end{displaymath}
Now the {\it prolongation} 
$\pro1\smfld\subset\Dx{n+1}$ of a submanifold 
$\smfld\subset\Dx n$ consists of the $n+1$ jets
$\jet{n+1}\mpnts\sigma$ with the property that the image
of the associated linear map  is tangent to $\smfld$, that is, $\lm_{\jet{n+1}\mpnts\sigma}(T_{\mpnts}M)\subset T_{\jet n\mpnts\sigma}\smfld$.\hbox{\vrule height 0pt depth 8pt width 0pt}

While it is customary to call a pseudogroup $\G\subset\D$ Lie if transformations $\locdif\in\G$ satisfy the condition, originally introduced by Lie \cite{LieI}, that they form the complete solution to a system of partial differential equations, several variants of the precise technical definition of a Lie pseudogroup exist in the literature, see e.g. \cite{GuilSter,HHJohnson, Kumpera,KuraI,SingSter}. For the purposes of this paper the following will suffice.
 
\begin{definition}\label{D:LiePG}
A subset $\G\subset\D$ is a {\it Lie pseudogroup}
if, whenever $\locdif$, $\psi\in\G$, then also
$\locdif\comp\psi^{-1}\in\G$ where defined, and, in addition, there is an integer $\ops\geq1$ so that for all $n\geq\ops$,
\begin{list}{}{\setlength{\leftmargin}{0.4in}\setlength{\parsep}{2pt}\setlength{\labelwidth}{0.15in}
}
\setlength{\hangindent}{30pt}
\item[\textbf{1.}] the corresponding subgroupoid 
$\Gx n\subset \Dx n$ forms a smooth, embedded subbundle;
\item[\textbf{2.}] every smooth function $\locdif\in\D$ satisfying $\jet n\mpnts\locdif\in\Gx n$, $\mpnt\in\dom\locdif$, belongs to $\G$;
\item[\textbf{3.}] $\Gx n = \pro{n - \ops}\Gx{\ops}$, $n\geq\ops$, agrees with the repeated prolongation of $\Gx{\ops}$.
\end{list}
\end{definition}

Thus on account of condition (1),  
for $n \geq \ops$, the pseudogroup subbundles 
$\Gx n\subset \Dx n$ are defined in local coordinates 
by formally integrable systems of $n^{\mbox{\scriptsize th}}$ 
order partial differential equations 
\begin{equation}\label{E:lPGDetEqns} 
F^{(n)}\zZn = 0,
\end{equation} 
the (local) \emph{determining equations} for the 
pseudogroup, whose local solutions $\mPnt = \locdif(\mpnt)$, 
by condition (2), are exactly the pseudogroup transformations.
Moreover, by condition (3), the determining equations in order 
$n>\ops$ can be obtained from those in order $\ops$ by a repeated 
application of the total derivative operators $\DD{\zc a}$ 
defined in \eqref{E:DTotDerOp}.

\emph{Remark}:  In \cite{IOV}, it is shown that, in the analytic category, the regularity condition (1) and Lie condition (2) imply the integrability condition (3).

Note that the customary requirements that a pseudogroup be closed under restriction of domains and concatenation of compatible local diffeomorphisms are built into condition (2). Thus our Lie pseudogroups are always complete in the sense of \cite{KuraII}. 
The assumptions also imply, as per the classical result 
of E.~Cartan \cite{Vara}, that the isotropy jets
\begin{equation}\label{E:Stabilizer}
\Gx{n}_{\mpnts} =
\{\gpx n\in\Gx n\;|\;\sourceDx n(\gpx n)=\targetDx n(\gpx n)=\mpnt\}\subset\ISD n{\mpnts}
\end{equation}
form a finite dimensional Lie group for all $\mpnt\in M$ and $n\geq\ops$.

Given a Lie pseudogroup $\G$, let $\g\subset\XM$ denote the set of its \emph{infinitesimal generators}, or $\G$ \emph{vector fields} for short. Thus $\g$ consists of the locally defined smooth vector fields $\vb$ on $M$ with the property that 
the flow maps $\flow{\vb}t$, for all fixed $t$, belong to $\G$.  
As a consequence of the group property in Definition \ref{D:LiePG}, the Lie bracket of two $\G$ vector fields, where defined, is again a $\G$ vector field.

Let $\gx n$ denote the space of the jets of $\G$ vector fields. In local coordinates \eqref{E:TMJetCoord}, the subspace 
$\gx n\subset \JTMx n$ 
 is specified by a linear 
system of partial differential equations  
\begin{equation}\label{E:LinDetEqns} 
\IDEn\zzcfn = 0,\qquad\jo\geq\ops,
\end{equation}
for the component functions $\vcx a{}=\vcx a{}(\zc{})$ of a vector field
obtained by linearizing the determining equations 
\eqref{E:lPGDetEqns} at the $n$-jet $\idjet n{\,\mpnts}=\jet n{\mpnts}\mkern1mu\id$ 
of the identity transformation.
Equations \eqref{E:LinDetEqns} are called the \emph{linearized} or 
\emph{infinitesimal determining equations} 
for the pseudogroup. As a consequence of Definition 
\ref{D:LiePG}, conversely, any vector field $\vb$ satisfying the infinitesimal 
determining equations \eqref{E:LinDetEqns}  can be shown to be 
an infinitesimal generator for $\G$, cf.~\cite{OPo}. Furthermore, as with the determining equations \eqref{E:lPGDetEqns} for pseudogroup transformations, the infinitesimal determining equations \eqref{E:LinDetEqns} in order $n\geq\ops$ can be obtained from those of order $\ops$ by repeated differentiation. 

While, by construction, the determining equations \eqref{E:lPGDetEqns} for a pseudogroup are locally solvable, that is, any $\zZn\in\Gx\jo$ is the jet of some $\locdiff\in\G$, it is not known to us if, in the $C^\infty$ category, the same holds true for the linearized version \eqref{E:LinDetEqns} of the equations. We will therefore make the additional blanket assumption that 
every $n$-jet $\vpntn\in\JTMx n$ satisfying \eqref{E:LinDetEqns} can be realized 
as the $n$-jet of some $\G$ vector field, that is, $\G$ is \emph{tame} as defined in \cite{OPo}. In this situation, the lift map
$\liftmapj\jo{\mpnts}$, as given in \eqref{E:LiftMapPoint}, restricts to an isomorphism between the Lie algebra 
$\gxfo n{\mpnts}$ of $n$-jets of $\G$ vector fields vanishing at $\mpnt$ and the 
Lie algebra of the isotropy subgroup $\displaystyle \IS{n}{\mpnts}$.

In the case of a symmetry group of a 
system of differential equations, the 
linearized determining equations \eqref{E:LinDetEqns} 
are the completion of the usual 
determining equations for the infinitesimal
symmetries obtained via Lie's algorithm 
\cite{OlverBook}. 

\section{Jet Bundles}\label{S:JB}

For $0\leq \jo\leq \infty $, let $\J\jo = \J\jo(M,\sdim)$ 
denote the $\jo^{\mbox{\scriptsize th}}$ order (extended) 
jet bundle consisting of equivalence 
classes of $\sdim$-dimensional submanifolds $S\subset M$ 
under the equivalence relation of $\jo^{\mbox{\scriptsize th}}$ 
order contact, cf.~\cite{Ehresmann, OlverBook}. We use the standard local coordinates 
\begin{equation}\label{E:zn}
\zx\jo = \xu\jo = (\xc i,\uc\alpha{},\uc\alpha{j_1},\uc\alpha{j_1j_2},\dots, \uc\alpha{j_1j_2\cdots j_\jo})
\end{equation}
on $\J\jo$ induced by a splitting of the local coordinates 
$\zc{} = (\xc{},\uc{}{})= (\sups {\xc{}{}} \sdim,\sups {\uc{}{}{}} \fdim)$
on $M=\J0$ into $\sdim$ independent and $\fdim=\mdim-\sdim$ dependent variables.
Let 
\[\projJ\gind\jo\colon\J \gind\longrightarrow\J n,\qquad 0\leq \jo\leq \gind,
\] 
denote the canonical projections. 

Local diffeomorphisms $\locdif\in\D$ preserve the 
$\jo^{\mbox{\scriptsize th}}$ order contact between submanifolds, and thus give rise to an action 
\begin{equation}\label{E:DactsJ}
\lactionJ{\locdif}(\jpx\jo)=\locdif\mkern -1mu\cdot\jpx\jo, \quad\mbox{where}\quad\jpx\jo\in(\projJ\jo o)^{-1}(\dom\locdif)\subset\J\jo,
\end{equation}
the so-called $\jo^{\mbox{\scriptsize th}}$ \emph{prolonged action}
of $\D$ on the jet bundle $\J\jo$. By the chain rule, the action \eqref{E:DactsJ} induces a well-defined action 
\begin{equation}\label{E:GactsJ}
\lactionJ{\gpx\jo}(\jpx\jo)=\gpx\jo\mkern -3mu\cdot\jpx\jo,\quad\mbox{where\quad
$\sourceDx\jo(\gpx\jo)=\projJ\jo o(\jpx\jo)$},
\end{equation}
of the diffeomorphism jet groupoid $\Dx n$ on $\J n$.

It will be useful to combine 
the two bundles $\Dx\jo$ and $\J\jo$ into a new bundle $\Ex\jo\to\J\jo$ 
by pulling back $\sourceDx\jo\colon\Dx\jo\to M$ via the standard projection 
$\pion\colon\J\jo\to M$. Thus $\Ex n$ consists of pairs of jets, 
\[
(\jpx\jo,\gpx \jo)\in\J\jo\times\Dx\jo,
\]
with $\jpx\jo\in\J\jo$ and $\gpx\jo\in\Gx\jo$ based at the same point $\mpnt=\projJ \jo o(\jpx\jo)=\sourceDx\jo(\gpx\jo)\in M$.

Local coordinates on $\Ex\jo$ are written as 
\begin{equation}\label{E:ELocalCoord}
\ZEn = \znZn,
\end{equation} 
where $\zx n = \xu n = (\xc i,\uc\alpha{},\uc\alpha{j_1},\uc\alpha{j_1j_2},\ldots,\uc\alpha{j_1j_2\cdots j_n})$ 
indicate submanifold jet coordinates, while 
\begin{displaymath}
\begin{split}
\Zn  &= (\Zc a{},\Zc a{b_1},\ldots,\Zc a{b_1b_2\cdots b_\jo})= (\Xn,\Un)\\
&=(\Xc i{},\Uc\alpha{},\Xc i{b_1},\Uc\alpha{b_1},\ldots,\Xc i{b_1b_2\cdots b_\jo},\Uc\alpha{b_1b_2\cdots b_\jo})
\end{split}
\end{displaymath} 
indicate the target derivative coordinates of a diffeomorphism. 
The source $\sourceEx n:\Ex n\to\J n$ and target $\targetEx n:\Ex n\to\J n$ maps on $\Ex n$ are respectively defined by
\begin{equation}\label{E:ESourceMap}
\sourceEx n\zngn = \jpx\jo,\qquad\targetEx\jo\zngn = \gpx\jo\mkern -3mu\cdot \jpx n.
\end{equation} 
Thus the latter simply represents the action of $\Dx n$ on $\J n$.

A local diffeomorphism $\locdif\in\D$ 
acts on the set 
\[
\{(\jpx n,\gpx n)\in\Ex n\,|\,\projJ n o(\jpx n) \in \dom\locdif\}\subset\Ex n
\] 
by
\begin{equation}\label{E:DActionE}
\lactionE\varphi\cdot\zngn =
(\jet n \mpnts\locdif\cdot\jpx n,
\gpx n\mkern -4mu\cdot \jet n{\locdif(\mpnts)}\locdif^{-1}),
\end{equation}
where $\projJ n o(\jpx n) = \mpnt$. The action
\eqref{E:DActionE} obviously factors into an action of $\Dx n$ on $\Ex n$, which we will again designate by the symbol $\lactionE{}$. Note that the target map $\targetEx n$ is manifestly invariant under the action  \eqref{E:DActionE} of the diffeomorphism pseudogroup, 
\begin{equation}\label{E:InvTargetMap}
\targetEx n(\lactionE\varphi\cdot\zngn)=\targetEx n\zngn.
\end{equation}

In local coordinates, the standard \emph{lifted total derivative operators} on $\Ei$ are given by
\begin{equation}\label{E:LDtot}
\DE{\xc j} = \DD {\xc j} + \sum _{\alpha =1}^\fdim 
\uc\alpha j \>\DD{\uc\alpha{}} + 
\sum_{\gind\,\geq\,1}\> \uc \alpha {jj_1j_2\cdots j_\gind}
\frac{\partial}{\partial \uc\alpha {j_1j_2\cdots j_\gind}}\,,
\end{equation}
where $\DD{\xc j}$, $\DD{\uc\alpha{}}$
are the total derivative operators 
\eqref{E:DTotDerOp} on $\Di$.
The \emph{lifted invariant total derivative operators} 
on $\Ei$ are, in turn, given by
\begin{equation}\label{E:DXi}
\DE{\Xc j{}} = \sum_{\gind=1}^p {\W \gind j} \, \DE{\xc\gind},
\qquad\where\qquad 
{\W \gind j} =(\DE{\xc\gind} {\Xc j{}})^{-1}
\end{equation}
indicates the entries in the inverse of the total 
Jacobian matrix; see \cite{OP2}. Then, by virtue of the chain rule, the expressions 
for the higher-order prolonged action of 
$\Dx n$ on $\J n$, that is, the coordinates 
$\Uhat\alpha J$ of the target map
$\targetEx n\colon \Ex n\to\J n$, are obtained by 
successively applying the derivative operators \eqref{E:DXi}
to the target dependent variables $\Uc\alpha{}$, 
\begin{equation}\label{E:UaJ}
\Uhat\alpha{j_1j_2\cdots j_\gind} 
= \DE{\Xc{j_1}{}}\DE{\Xc{j_2}{}}\cdots  \DE{\Xc{j_\gind}{}} \Uc\alpha{}.
\end{equation}
Note that we employ hats in \eqref{E:UaJ} to distinguish between the target jet coordinates of submanifolds and diffeomorphisms. 

Let $\vb\subset\X(M)$ be a smooth vector field 
with the flow map $\flow {\vb} t$. By definition,
the prolongation $\prx\jo\vb$ of $\vb$ is the infinitesimal generator 
of the prolonged action of $\flow{\vb}t$ on 
$(\projJ\jo o{})^{-1}(\dom\vb)\subset \J n$. 
Write 
\[
\vb=\sum_{i=1}^\sdim\vxx i \dfrac{\partial}{\partial \xc i}+
\sum_{\alpha=1}^{\fdim}\vux\alpha\dfrac{\partial}{\partial\uc\alpha{}}
\]
in the coordinates \eqref{E:zn}. Then
the components $\prc \alpha {j_1j_2\cdots j_\gind}$ of
\begin{equation}\label{eq:prdefn}
\prx\jo\vb = \sum_{i=1}^\sdim\vxx i \dfrac{\partial}{\partial \xc i}+
\sum_{\alpha=1}^{\fdim}\sum_{\gind\leq\jo}\prc\alpha{j_1j_2\cdots j_\gind}\dfrac{\partial}{\partial\uc\alpha{j_1j_2\cdots j_\gind}}
\end{equation}
are given by the standard prolongation formula 
\begin{equation}\label{E:prncomp}
\prc\alpha{j_1j_2\cdots j_\gind} =  \DEder{\xc{j_1}}\DEder{\xc{j_2}}\cdots\DEder{\xc{j_\gind}}\charv\alpha+\sum_{i=1}^{\sdim}\vxx i\uc\alpha{ij_1j_2\cdots j_\gind},\qquad
\end{equation}
where 
\begin{equation}\label{eq:char}
\charv\alpha = \vux\alpha{}-\vxx i\uc\alpha i,\qquad
\alpha = 1,\dots,\fdim,
\end{equation}
denotes the components of the characteristic of $\vb$ and $\DE{\xc j}$ stands for the total derivative operators \eqref{E:LDtot} restricted to $\Ji$, identified as the image of the identity section 
\[
\Ex{\infty}{}_{|\idjet{\infty}{}}=\{(\jpx\infty,\idjet \infty\mpnts)\;\vert\;\jpx\infty\in\J\infty,\mpnt=\projJ\infty o(\jpx\infty)\}
\]
in $\Ex\infty$, cf. \cite{OlverBook}.

Finally, in view
of \eqref{E:prncomp}, the 
prolongation $\prx{n}\vb(\jpx\jo)$
of a vector field
at $\jpx\jo\in\J n$ depends only on the $\jo$-jet $\jet\jo\mpnts\vb$ of
$\vb$ at $\mpnt=\projJ\jo o(\jpx\jo)$ and, consequently,
the prolongation process induces well-defined linear mappings
\begin{equation}\label{E:Prolongation}
\prm_{\jpxs n}: \JTMxf n\mpnts\>\longrightarrow\>
T_{\jpxs\jo}\J\jo,\qquad \mpnts = \projJ n o(\jpx\jo). 
\end{equation}

\section{Moving Frames}\label{S:MF}

Given a Lie pseudogroup $\G\subset\D$, we let $\Hx\jo\subset\Ex\jo$ denote the subbundle corresponding to the jets of transformations belonging to $\G$. Specifically,
\begin{equation}
\Hx\jo = \{(\jpx\jo,\gpx\jo)\in\Ex\jo\,\vert\,\gpx\jo\in\Gx\jo\}.
\end{equation}
We will furthermore designate the restrictions of the source and target maps \eqref{E:ESourceMap} to $\Hx n$ by $\sourceHx n$, $\targetHx n$.
Let $\Uset\subset\J n$ be open and connected. Then
a \emph{local moving frame} $\mfx n$ on $\Uset$ for the action of $\G$ on $\J\jo$ 
is a section of 
\[
\sourceHx n\colon\Hx n{}\mkern-2mu{_{|\Uset}}\>\longrightarrow\> \Uset
\]
that is locally equivariant, that is, there is an open set
\begin{equation}\label{E:Wset}
\Wset\subset(\sourceHx n)^{-1}(\Uset)\cap(\targetHx n)^{-1}(\Uset)
\end{equation}
containing the image of the identity section 
$\{(\jpx n,\idjet n{\mpnts})\;\vert\;\jpx n \in\Uset\}\subset\Wset$
so that 
\begin{equation}\label{E:Equivariance}
\mfx n(\gpx n\mkern-3mu\cdot\jpx n) = \lactionE{\gpx n}\rho(\jpx n),
\qquad\mbox{for all $(\jpx n,\gpx n)\in\Wset$}.
\end{equation}
Note that if \eqref{E:Equivariance} holds in the open sets
$\Wset_1$, $\Wset_2\subset\Hx n$, then it also holds in the union
$\Wset_1\cup\Wset_2$, so that one can always assume
that $\Wset$ is the maximal set with the required
properties. 

A section of $\Hx n{}_{|\Uset}\to\Uset$ is called a {\it global moving frame}, or simply a {moving frame}, if $\Uset$ is stable under the action of $\Gx n$, that is, $\Uset$ is the union of the orbits of the $\Gx n$ action on $\J n$, and if $\Wset$ in the equivariance condition
\eqref{E:Equivariance} can be chosen to be the entire set $\Wset = \Hx n{}\mkern-2mu{_{|\Uset}}$.
Call a local moving frame $\mfx n\colon\Uset\to\Hx n$ \emph{normalized} if 
$\mfx n(\jpx n) = (\jpx n,\idjet n\mpnts)$ for some $\jpx n\in\Uset$.

Moving frames $\mfx n_1\colon\Usetx n\to\Hx n$, 
$\mfx k_2\colon\Usetx k\to\Hx k$, $k>n$, 
are said to be \emph{compatible} if $\projJ k n(\Usetx k)=\Usetx n$ and   
\begin{equation}\label{E:compMF}
\mfx n_1\comp\projJ k n(\jpx k)=\projE k n\comp\mfx k_2(\jpx k)
\end{equation}
for all $\jpx k\in\Usetx k$, where $\projE k n\colon\Hx k\to \Hx n$ stands for the canonical projection.
A \emph{complete moving frame} is provided 
by the projective limit of a mutually compatible collection 
$\mfx k\colon\Usetx k\to \Hx k$ of moving frames 
of all orders $k\geq n$ for some $n$.
As expounded in \cite{OP2}, complete moving frames can be effectively used 
to construct complete sets of differential invariants, 
invariant total derivative operators, invariant coframes, and so on, and to analyze the structure of the algebra of differential invariants
for the action of pseudogroups on extended jet bundles.   

As for Lie transformation groups \cite{FelsOlver1999}, the existence of a moving frame hinges on a suitable notion of freeness of the pseudogroup action on the jet bundles $\J\jo$. However, in contrast with the finite dimensional case, bona fide infinite dimensional transformation groups cannot have trivial isotropy, and, as a result, we are lead to define freeness of the action in terms of jets of local diffeomorphisms stabilizing a given submanifold jet.

Recall that as a consequence of Definition \ref{D:LiePG}, the \emph{isotropy subgroup} 
\[\IS\jo {\jpxs\jo}= \{\gpx\jo\in \IS{n}{\mpnts}\;|\;
\gpx\jo\mkern-2mu\cdot\mkern-2mu \jpx\jo = \jpx\jo\}\] of a point $\jpx\jo\in\J n$, as a closed subgroup,
forms a Lie subgroup of $\displaystyle\IS{\jo}{\mpnts}$, where $\mpnt=\projJ\jo o(\jpx\jo)$. In addition, one can show that the Lie algebra of 
$\IS \jo{\jpxs\jo}$ can be identified with the kernel of the restriction of the 
prolongation map $\prm_{\jpxs\jo}$ in \eqref{E:Prolongation} to $\gxfo \jo\mpnts$.

\begin{definition}\label{D:freedf}
A pseudogroup $\G$ acts \emph{freely} at 
$\displaystyle\jpxo\jo\in \J\jo$ 
if its  \emph{isotropy subgroup} is trivial,
$\IS n{\jpxos n}= 
\{\idjet n\mpntos\}$,
and \emph{locally freely}  
if $\IS n{\jpxos n}$ is discrete. 
\end{definition}
Thus the pseudogroup $\G$ acts locally freely at $\displaystyle\jpxo\jo$ precisely when the prolongation map $\prm_{\jpxos n}\colon\gxf n\mpntos\to T_{\jpxos n}\J n$ is injective. In this situation the mappings $\prm_{\jpxs\jo}\colon\gxf \jo\mpnts\to T_{\jpxs\jo}\J\jo$ have maximal rank for all $\jpx\jo$ contained in some neighborhood $\VsetT{}\subset\J\jo$ of  $\displaystyle\jpxo\jo$ and thus their images define an involutive distribution on $\VsetT{}$ whose integral submanifolds are the intersections of $\G$-orbits on $\J\jo$ with $\VsetT{}$. A  \emph{cross section}, or \emph{transversal}, $\crosss\jo$ to the orbits of $\G$ through $\displaystyle\jpxo\jo$ is an embedded submanifold of $\J n$ containing $\displaystyle\jpxo n$ so that 
\begin{equation}\label{E:Transversal}
T_{\jpxs\jo}\J n= T_{\jpxs\jo}\crosss\jo\oplus \range\prm_{\jpxs \jo},\qquad \mbox{for all $\jpx\jo\in\crosss\jo$}.
\end{equation}
Note that the existence of cross sections for locally free actions is a simple consequence of the classical Frobenius theorem, \cite{OlverBook}.

\begin{theorem}\label{T:mf} 
Suppose $\G$ acts locally freely at 
$\displaystyle\jpxo\jo\in\J \jo$. 
Then $\G$ admits a normalized local moving frame on some neighborhood $\Uset\subset\J\jo$ of 
$\displaystyle\jpxo \jo$. 
Suppose furthermore 
that one can choose a cross section $\crosss\jo$ through 
$\displaystyle\jpxo \jo$ 
so that $\Gx\jo$ acts freely at each 
$\cpx\jo\in\crosss\jo$ 
and that any $\G$-orbit intersects $\crosss\jo$ in at most one point. Then $\G$ admits a global moving frame in some open set $\Uset\subset\J\jo$ containing $\displaystyle\jpxo\jo$.
\end{theorem} 
 
\begin{proof} By assumption, the mappings $\prm_{\jpxs\jo}\colon\gxf \jo\mpnts\to T_{\jpxs\jo}\J\jo$ have maximal rank for all $\jpx\jo$ contained in some neighborhood $\VsetT{}\subset\J\jo$ of  $\displaystyle\jpxo\jo$. 
Let $\crosss\jo\subset\VsetT{}$ be a cross section to the orbits through $\displaystyle\jpxo\jo$
and write $\Hx\jo{}_{|\crosss\jo} = (\sourceHx\jo)^{-1}(\crosss\jo)$. Let
\begin{equation}\label{E:targetcross}
\targetHCx\jo = \mathstrut{\targetHx\jo{}}_{|\crosss\jo}\colon\Hx\jo{}_{|\crosss\jo}\>\longrightarrow\>\J \jo
\end{equation}
denote the target map restricted to $\Hx\jo{}_{|\crosss\jo}$.
By \eqref{E:Transversal}, the Jacobian of $\targetHCx\jo$ is non-singular at $\displaystyle(\jpxo\jo,\idjet\jo\mpntos)$, and so, by the inverse function theorem, $\targetHCx\jo$ restricts to a diffeomorphism from a neighborhood $\Vset{}\subset\Hx\jo{}_{|\crosss\jo}$ of $\displaystyle(\jpxo\jo,\idjet\jo\mpntos)$ onto a neighborhood $\Uset\subset\J\jo$ of $\displaystyle\jpxo\jo$. Write 
$\targetHCinvx\jo=(\iota^{(\jo)},\gamma^{(\jo)})\colon\Uset\to\Vset{}$ for the inverse function and define a section $\mfx\jo\colon\Uset\to\Hx\jo$ by
\[
\mfx\jo(\jpx\jo) = (\jpx\jo,\gamma^{(\jo)}(\jpx\jo)^{-1}),
\]
where the exponent indicates the groupoid inverse on $\Gx\jo$. A direct computation shows that for $\jpx n=\targetHCx\jo(\cpx\jo,\hpx \jo)$, where $(\cpx\jo,\hpx \jo)\in\Vset{}$, the equivariance condition
\begin{equation}\label{E:rhoEquivariance}
\mfx\jo(\gpx\jo\mkern-3mu\cdot
\jpx\jo) = \lactionE{\gpxs\jo}\rho(\jpx\jo)
\end{equation}
is satisfied provided that $(\cpx\jo,\gpx\jo\mkern-3mu\cdot\hpx\jo)\in\Vset{}$.
But it is easy to see that the pairs $(\jpx\jo,\gpx\jo)$ fulfilling this condition form an open set $\Wset\subset\Hx n{}_{|\Uset}$ containing the image of the identity section, and, consequently, $\mfx\jo$ provides a normalized local moving frame in the neighborhood $\Uset$ of $\displaystyle\jpxo\jo$.

Next assume that $\crosss\jo$ is a cross section through $\displaystyle\jpxo n$ so that $\Gx\jo$ acts freely at every $\cpx\jo\in\crosss\jo$ and that each $\Gx\jo$ orbit intersects $\crosss\jo$ at most at one point. These conditions are equivalent to the mapping $\targetHCx\jo$ defined in \eqref{E:targetcross} being one-to-one, and so the steps used above to construct a local moving frame will also prove the existence of the global counterpart provided that the rank of $\targetHCx\jo$ is maximal at every point.

To compute the rank of $\targetHCx\jo$ at $(\cpxo\jo,\hpxo\jo)\in\Hx\jo{}_{|\crosss\jo}$, write $\hpxo\jo=\jet\jo\mpntos\locdif$, $\locdif\in\G$, and consider the mapping
\begin{equation}\label{eqn:emult}
\begin{split}
\emult\locdif:&(\projJ\jo o\comp\targetHCx\jo)^{-1}(\dom\locdif)\subset\Hx\jo{}_{|\crosss\jo}\longrightarrow\>\Hx\jo{}_{|\crosss\jo};\\ &\qquad\emult\locdif(\cpx\jo,\hpx\jo) = (\cpx\jo,\jet\jo{\targetDx\jo (\hpxs\jo)}\locdif\cdot\hpx\jo).
\end{split}
\end{equation}
Then obviously
\[
\targetHCx\jo\comp\emult\locdif=\lactionJ\locdif\comp\targetHCx\jo,
\]
so
\[
\targetHCx\jo_*{}\mid_{(\cpxos\jo,\hpxos\jo)}\comp\emult\locdif_*{}\mid_{(\cpxos\jo,\idjet n{\mpntos})} =(\lactionJ\locdif)_*{}\mid_{\cpxos\jo}\comp\targetHCx\jo_*{}\mid_{(\cpxos\jo,\idjet n{\mpntos})}.
\]
By assumption, the ranks of the differentials on the right-hand side of the equation are maximal, so $\targetHCx\jo$ must indeed have maximal rank at $(\cpxo\jo,\hpxo\jo)\in\Hx\jo{}_{|\crosss\jo}$. This completes the proof of the theorem. 
\end{proof}

\section{Persistence of Freeness}\label{S:PF}

In this final section we state and prove the main results of the paper establishing the persistence of both local and global freeness under prolongation of the pseudogroup action. 

\begin{theorem}\label{T:LocalFreeness} Suppose $\G$ acts locally freely at $\displaystyle\jpxo\jo$, where $\jo\geq\ops$. Then it acts locally freely at any $\displaystyle\jpxo{\jo+\gind}\in\J{\jo+\gind}$ with $\displaystyle\projJ{\jo+\gind}{\jo}(\jpxo{\jo+\gind}) = \jpxo\jo$. 
\end{theorem}

\begin{proof}
It suffices to consider the case $k=1$ only. Fix $\displaystyle\jpxo{\jo+1}\in\J{\jo+1}$ with $\displaystyle\projJ{\jo+1}{\jo}(\jpxo{\jo+1}) = \jpxo\jo$. Recall that $\G$ acts locally freely at $\displaystyle\jpxo{\jo+1}\in\J{\jo+1}$ if and only if the restriction of the prolongation map 
$$\prm_{\jpxos{\jo+1}}\colon\gxf{\jo+1}\mpntos\>\longrightarrow\> T_{\jpxos{\jo+1}}\J{\jo+1}$$
 is injective, that is, 
\begin{equation}\label{E:LocFreeCond}
\gxf {\jo+1}\mpntos\cap\ker\prm_{\jpxos{\jo+1}}=\{0\}.
\end{equation}

Let $\displaystyle\vpntxo{\jo+1}\in\gxf {\jo+1}\mpntos\cap\ker\prm_{\jpxos{\jo+1}}$. 
Then obviously the projection $\displaystyle\vpntxo{\jo}$ of $\displaystyle\vpntxo{\jo+1}$ into $\JTMx\jo$ satisfies
\begin{displaymath}
\vpntxo{\jo}\in\gxf {\jo}\mpntos\cap\ker\prm_{\jpxos{\jo}},
\end{displaymath}
so, by assumption, $\displaystyle\vpntxo{\jo}$ must vanish. Thus in local coordinates,

\begin{equation}\label{E:vCoord}
\vpntxo{\jo+1}=(\zoc a, 0,\dots,0,\zcfxo b{c_1c_2\cdots c_{\jo+1}}),  
\end{equation} 
where the components $\zcfxo b{c_1c_2\cdots c_{n+1}}$ are determined by the requirements that the jet $\displaystyle\vpntxo{n+1}$ satisfy the infinitesimal determining equations \eqref{E:LinDetEqns} of order $n+1$ and be contained in the kernel of the prolongation map  $\ker\prm_{\jpxos{n+1}}$. 

Recall that the infinitesimal determining equations in order $n+1$ can be obtained from those in order $n$ by differentiation. Thus an equation 
\[
\sum_{0\leq\gind\leq\jo}\sum_{1\leq b,c_1,\ldots,c_\gind\leq \mdim}\IDE{\lower1.5pt\hbox{$\scriptstyle A$}} b{c_1c_2\cdots c_\gind}(\zc a)\zcfx b{c_1c_2\cdots c_\gind}=0
\]
of order $\jo$ yields the equations 
\begin{equation}\label{E:LinDetEqnsSym}
\sum_{1\leq b,c_1,\ldots,c_\jo\leq \mdim}\IDE{\lower1.5pt\hbox{$\scriptstyle A$}} b{c_1c_2\cdots c_\jo}(\zoc a)\zcfxo b{c_1c_2\cdots c_\jo c_{\jo+1}}=0,\quad c_{\jo+1}=1,\dots,\mdim,
\end{equation}
for the coordinates $\zcfxo b{c_1c_2\cdots c_{\jo+1}}$, and $\displaystyle\vpntxo{\jo+1}\in\gx{\jo+1}$ precisely when all the derived equations of this form are satisfied.  

Divide, as usual, the local coordinates $(\zc{a}) =(\xc{i},\uc{\alpha}{})$ of $M$ into independent and dependent variables, and denote the induced coordinates on $\JTMxf\jo{}$ by
\begin{displaymath}
(\zcfx a{c_1c_2\cdots c_\gind}) = (\zxc i{c_1c_2\cdots c_\gind},\zuc \alpha{c_1c_2\cdots c_\gind}).
\end{displaymath}
Next define the differential operators 
\begin{equation}\label{E:PrSymb}
\dop i = \DT{\xc i}+\uco\alpha i\DT{\uc \alpha{}}
\end{equation}
on $\JTMxf\infty{}$, where $\DT{\xc i}$, $\DT{\uc \alpha{}}$ denote the standard total derivative operators on $\JTMxf\infty{}$ and $\uco\alpha{i}$ is the (constant) first order derivative coordinate of the jet $\displaystyle\jpxo{\jo+1}$. 

Then, on account of \eqref{E:vCoord}, the components of the derivative variables in $\displaystyle\prm_{\jpxos{n+1}}\vpntxo{n+1}$ of order $k\leq n$ vanish, while the vanishing of the $\uc\alpha{i_1i_2\cdots i_{n+1}}$-component $\prc\alpha{o,i_1i_1\cdots i_{n+1}}$ of $\displaystyle\prm_{\jpxos{n+1}}\vpntxo{n+1}$ yields the equations 
\begin{equation}\label{E:KerPREqns}
\prc\alpha{o,i_1i_2\cdots i_{n+1}} = \dop{i_1}\dop{i_2}\cdots\dop{i_{n+1}}
(\prvx\alpha{}-\uco\alpha j\vxx j)=0
\end{equation}
for the coordinates $\zcfxo b{c_1c_2\cdots c_{n+1}}$ of $\displaystyle\vpntxo{n+1}$.

Fix $1\leq i\leq\sdim$, and let $\wpntx \jo_{i}\in\JTMxf \jo\mpntos$ denote the jet with coordinates
\begin{equation}\label{E:wpoint}
\wpntx \jo_{i} = (\zoc a,0,\dots,0,\zcfx b{c_1c_2\cdots c_\jo}=\zcfxo b{ic_1c_2\cdots c_\jo}+\uco\alpha i\zcfxo b{\alpha c_1c_2\cdots c_\jo}).
\end{equation}
Then, on account of equations \eqref{E:LinDetEqnsSym}, \eqref{E:KerPREqns}, we have that 
\[
\wpntx\jo_{i}\in\gxf {\jo}\mpntos\>\cap\>\ker\prm_{\jpxos{\jo}}=\{0\},
\]
 and so, by the assumptions,
\begin{equation}\label{E:dadaa}
\zcfxo b{ic_1c_2\cdots c_\jo}+\uco\alpha i\zcfxo b{\alpha c_1c_2\cdots c_\jo}=0,\qquad i=1,\ldots,\sdim.
\end{equation}

Finally, let $\wpntxhat\jo_{e}\in\JTMxf \jo\mpntos$, $1\leq e\leq\mdim$, be the jet with coordinates
\begin{equation}
\wpntxhat \jo_{e} = (\zoc a,0,\dots,0,\zcfxo b{c_1c_2\cdots c_\jo e}).
\end{equation}
Then, by virtue of \eqref{E:LinDetEqnsSym}, \eqref{E:dadaa}, we have that 
\[
\wpntxhat\jo_e\in\gxf{\jo}\mpntos\>\cap\>\ker\prm_{\jpxos{\jo}}=\{0\}.
\] 
Consequently, $\zcfxo b{c_1c_2\cdots c_{\jo}e}=0$, which concludes the proof of the Theorem.
\end{proof}
\begin{theorem}\label{T:GlobalFreeness} 
Suppose $\G$ acts freely at $\displaystyle\jpxo n$, where $\jo\geq \ops+1$. Then it acts freely at any $\displaystyle\jpxo{\jo+\gind}\in\J{\jo+\gind}$ with $\displaystyle\projJ{\jo+\gind}{\jo}(\jpxo{\jo+\gind}) = \jpxo n$. 
\end{theorem}

\begin{proof} It suffices to prove that $\G$ acts freely at any $\displaystyle\jpxo{\jo+1}\in\J{\jo+1}$ with $\displaystyle\projJ{\jo+1}{n}(\jpxo{\jo+1}) = \jpxo n$.
Let $\displaystyle\gpxo{\jo+1}\in \IS{\jo+1}{\jpxos{\jo+1}}$. 

Then obviously 
$\displaystyle\projD{\jo+1}n(\gpxo{\jo+1})\in\IS{\jo}{\jpxos n}$, so
by assumption, $\displaystyle\gpxo{\jo+1}$ agrees with the jet of the identity transformation 
up to order $\jo$, that is, $\displaystyle\projD{\jo+1}n(\gpxo{\jo+1})=\idjet \jo{\mpntos{}}$. Thus in local coordinates,
\begin{equation}\label{E:goCoord}
\begin{split}
\gpxo{n+1}= (\zc a=\zoc a, \Zc a{}=\zoc a, \Zc a b &= \delta^a_b,\Zc a{b_1b_2}=0,\dots,\\
\Zc a{b_1\cdots b_\jo}&=0,\Zc a {b_1\cdots b_{\jo+1}}=\Zoc a{b_1\cdots b_{\jo+1}})
\end{split}
\end{equation}
for some $\Zoc a{b_1\cdots b_{\jo+1}}$. These
coordinates are determined by two sets of equations,
the first one specifying that $\displaystyle\gpxo{n+1}$ belongs to
$\Gx{n+1}$ and the second one imposing the condition
that the transformation on the fiber $\Jxf{n+1}{\mpntos{}{}}$ induced by $\displaystyle\gpxo{n+1}$ fixes $\displaystyle\jpxo{n+1}$.

We start with the first set of conditions. Since $n\geq\ops+1$, we can, on account of Definition \ref{D:LiePG}, prolong the determining equations for $\G$ in order $n-1$ to conclude that 
there is a neighborhood $\Vset {}\subset\Dx n$ of $\idjet n{\mpntos}$ so that 
$\Gx n\cap\Vset{}$ is the solution set of a system
of equations of the form
\begin{equation}\label{E:DetEqns}
\begin{split}
F_{\lower1.2pt\hbox{$\scriptstyle \alpha$}}{}_{,a}^{\hphantom{,}b_1\cdots b_n}\zZc{n-1}\Zc a{b_1\cdots b_n} + G_\alpha\zZc{n-1}&=0,\\
H_\beta\zZc{n-1}&=0.
\end{split}
\end{equation}
Furthermore, condition 3 of Definition \ref{D:LiePG} stipulates that, in addition to system \eqref{E:DetEqns},
pseudogroup jets $\gpx{n+1}\in\Gx {n+1}$ are determined by the equations
\begin{equation}\label{E:DetEqnsProl}
\begin{split}
F_{\lower1.1pt\hbox{$\scriptstyle \alpha$}}{}_{,a}^{\hphantom{,}b_1\cdots b_n}\zZc{n-1}&\Zc a{b_1\cdots b_ne}
+(\DD e F_\alpha,{}_a^{b_1\cdots b_n})\zZn\Zc a {b_1\cdots b_n}\\
&\quad+(\DD e G_\alpha)\zZn=0,\qquad e=1,\ldots,m,
\end{split}
\end{equation}
on the entire cylinder $\Vset1=(\projD{n+1}n)^{-1}(\Vset{})$. 

Now evaluate equations \eqref{E:DetEqnsProl} at $\displaystyle\gpxo{n+1}$ as given in \eqref{E:goCoord} to see that
\begin{equation}\label{E:DetEqnsProlSubs}
\begin{split}
&F_{\lower1.2pt\hbox{$\scriptstyle \alpha$}}{}_{,a}^{\hphantom{,}b_1\cdots b_n}(\zoc c,\zoc c,\delta^c_d,0,\dots,0)\Zoc a{b_1\cdots b_ne}\\
&\qquad+\dfrac{\partial G_\alpha}{\partial \zc e}(\zoc c,\zoc c,\delta^c_d,0,\dots,0)
+\dfrac{\partial G_\alpha}{\partial \Zc e{}}(\zoc c,\zoc c,\delta^c_d,0,\dots,0)=0.
\end{split}
\end{equation}
On the other hand, equation \eqref{E:DetEqns}, when evaluated at the identity jet $\idjet n\mpnts$ becomes
\begin{equation}
G_\alpha(\zc c,\zc c,\delta^c_d,0,\dots,0)=0,
\end{equation}
which, after differentiation with respect to $\zc e$, shows that
\eqref{E:DetEqnsProlSubs} reduces to a system of linear, homogeneous  
equations 
\begin{equation}\label{E:HODetEqns}
F_{\lower1.2pt\hbox{$\scriptstyle \alpha$}}{}_{,a}^{\hphantom{,}b_1\cdots b_n}(\zoc c,\zoc c,\delta^c_d,0,\dots,0)\Zoc a{b_1\cdots b_ne}=0,\qquad e=1,\dots,m,
\end{equation}
for the coordinates $\Zoc a{b_1\cdots b_{n+1}}$.

Next we use formulas \eqref{E:UaJ} to compute
the action of $\displaystyle\gpxo{n+1}$ at $\displaystyle\jpxo{n+1}$. The
components of interest are those in order $n+1$,
and these are given by
\begin{equation}\label{E:ProlAction}
\begin{split}
\Uhat\alpha{j_1\cdots j_{n+1}} &=
\DEder{X^{j_1}}\cdots\DEder{X^{j_{n+1}}} \Uc\alpha{}\\
&=(\W{j_1}{k_1}\DEder{x^{k_1}})\cdots 
(\W{j_{n+1}}{k_{n+1}}\DEder{x^{k_{n+1}}})\Uc\alpha{}.
\end{split}
\end{equation}
On account of \eqref{E:goCoord}, the only non-zero terms in \eqref{E:ProlAction} arise from 
\begin{equation}
\begin{split}
&\W{j_1}{k_1}\cdots \W{j_{n+1}}{k_{n+1}}\DE{x^{k_1}}
\cdots\DE{x^{k_{n+1}}}U^\alpha\qquad\text{and}\\
&\W{j_1}{k_1}\cdots \W{j_{n}}{k_{n}} (\DE{x^{k_1}}\cdots 
\DE{x^{k_{n}}}W_{j_{n+1}}^{k_{n+1}})
(\DE{x^{k_{n+1}}}U^\alpha).
\end{split}
\end{equation}
After some manipulations we see that
\begin{equation}\label{E:LinStabEqns}
\begin{split}
&\Uhat\alpha{j_1\cdots j_{n+1}}\\& =
\uco\alpha{j_1\cdots j_{n+1}}+
(\DD{x^{j_1}}+\uco{\gamma_1}{j_1}\DD{\uc{\gamma_1}{}})\cdots (\DD{x^{j_{n+1}}}+\uco{\gamma_{n+1}}{j_{n+1}}\DD{\uc{\gamma_{n+1}}{}})\Uc\alpha{}\hskip-10pt\\
&-\uco\alpha{k_{n+1}}(\DD{x^{j_1}}+\uco{\gamma_1}{j_1}
\DD{\uc{\gamma_1}{}})\cdots (\DD{x^{j_{n+1}}}+\uco{\gamma_{n+1}}{j_{n+1}}\DD{\uc{\gamma_{n+1}}{}})X^{k_{n+1}},\hskip-10pt
\end{split}
\end{equation}
where $\uco\alpha{j_1\cdots j_\gind}$ denote the coordinates of $\displaystyle\jpxo{n+1}$. 
Hence the conditions $\Uhat\alpha{j_1\cdots j_{n+1}} =\uco\alpha{j_1\cdots j_{n+1}}$ lead to another
system of linear, homogeneous equations for the coordinates
$\Zoc a{b_1\cdots b_{n+1}}$ in addition to \eqref{E:HODetEqns}.

Since $\G$ acts freely at $\displaystyle\jpxo{\jo}$, it also acts locally freely at $\displaystyle\jpxo{\jo}$, and, consequently, also at $\displaystyle\jpxo{n+1}$ by Theorem \ref{T:LocalFreeness}. This implies that the solution set to the homogeneous linear system obtained by combining \eqref{E:HODetEqns} and the equations resulting from \eqref{E:LinStabEqns} must discrete.  Consequently, it must be trivial, and hence $\G$ acts freely at $\displaystyle\jpxo{n+1}$.
\end{proof}


\begin{thebibliography}{99}
\footnotesize\itemsep=0pt

\bibitem{Cartan1937} Cartan, \'E., \emph{La Th\'eorie des Groupes Finis et Continus et la G\'eom\'etrie Diff\'erentielle Trait\'ees par la M\'ethode du Rep\`ere Mobile}, Gauthier-Villars, Paris, 1937.

\bibitem{CartanI} Cartan, \'E., \emph{Sur la structure des groupes infinis de transformations}, Oeuvres Compl\`etes, Part. II, Vol. 2, Gauthier--Villars, Paris, 1953, pp. 571--714.

\bibitem{COP1} Cheh, J., Olver, P.J., Pohjanpelto, J., \emph{Maurer--Cartan equations for Lie symmetry pseudo-groups of differential equations}, J. Math. Phys. 46 (2005) 023504.

\bibitem{COP2} Cheh, J., Olver, P.J., Pohjanpelto, J., \emph{Algorithms for differential invariants of symmetry groups of differential equations}, Found. Comput. Math. 8 (2008), 501-–532.\par


\bibitem{Ehresmann} Ehresmann, C., \emph{Introduction \`a la th\'eorie des structures infinit\'esimales et des pseudo-groupes de Lie}, G\'eometrie Diff\'erentielle; Colloq. Inter. du Centre Nat. de la Rech. Sci., Strasbourg (1953), 97--110.

\bibitem{FelsOlver1999} Fels, M., Olver, P.J., \emph{Moving coframes. II. Regularization and theoretical foundations}, Acta Appl. Math. 55 (1999), 127--208.

\bibitem{Green1978} Green, M.L., \emph{The moving frame, differential invariants and rigidity theorems for curves in homogeneous spaces}, Duke Math. J. 45 (1978), 735--779.

\bibitem{Griffiths1974} Griffiths, P., \emph{On Cartan's method of Lie groups and moving frames as applied to existence and uniqueness questions in differential geometry}, Duke Math. J. 41 (1974), 775--814.

\bibitem{GuilSter} Guillemin, V., Sternberg, S., \emph{Deformation Theory of Pseudogroup Structures}, Mem. Amer. Math. Soc. 64, AMS, Providence, RI, 1966.

\bibitem{IOV} Itskov, V., Olver, P.J., Valiquette, F., \emph{Lie completion of pseudo-groups}, preprint, University of Minnesota, 2009.

\bibitem{Jensen1977} Jensen, G.R., \emph{Higher Order Contact of Submanifolds of Homogeneous Spaces}, Lecture
Notes in Math., No. 610, Springer–Verlag, New York, 1977.

\bibitem{HHJohnson} Johnson, H.H., \emph{Classical differential invariants and applications to partial differential equations}, Math. Ann. 148 (1962), 308--329.

\bibitem{KOivb} Kogan, I.A., Olver, P.J., \emph{Invariant Euler-Lagrange equations and the invariant variational bicomplex}, Acta Appl. Math. 76 (2003), 137--193.

\bibitem{Kumpera} Kumpera, A., \emph{Invariants diff\'erentiels d'un pseudogroupe de Lie}, J. Diff. Geom. 10 (1975), 289--416.

\bibitem{Mackenzie} Mackenzie, K., \emph{General Theory of Lie Groupoids and Lie Algebroids}, London Math. Soc. Lecture Notes, vol. 213, Cambridge Univ.~Press, 
Cambridge, 2005.

\bibitem{KuraI} Kuranishi, M., \emph{On the local theory of continuous infinite pseudo groups I}, Nagoya Math. J. 15 (1959), 225--260.

\bibitem{KuraII} Kuranishi, M., \emph{On the local theory of continuous infinite pseudo groups II}, Nagoya Math. J. 19 (1961), 55--91.

\bibitem{LieI} Lie, S., \emph{\"Uber unendlichen kontinuierliche Gruppen}, Christ. Forh. Aar. 8 (1883), 1--47; also {\it Gesammelte Abhandlungen}, Vol. 5, B.G. Teubner, Leipzig, 1924, {pp.~314--360}.


\bibitem{Moro2005} Morozov, O.I., \emph{Structure of symmetry groups via Cartan's method: survey of four approaches}, SIGMA 1 (2005), 006.

\bibitem{OlverBook} Olver, P.J., \emph{Applications of Lie Groups to Differential Equations}, 2nd ed., Graduate Texts in Mathematics, vol. 107, Springer--Verlag, New York, 1993.

\bibitem{OlverEIS} Olver, P.J., \emph{Equivalence, Invariants, and Symmetry}, Cambridge University Press, Cambridge, 1995.

\bibitem{Olversurvey} Olver, P.J., \emph{A survey of moving frames}, in: \emph{Computer Algebra and Geometric Algebra with Applications}, H. Li, P.J. Olver and G. Sommer, eds., Lecture Notes in Computer Science, vol. 3519, Springer--Verlag, New York, 2005, pp. 105--138.

\bibitem{OPo} Olver, P.J., Pohjanpelto, J., \emph{Regularity of pseudogroup orbits}, in Symmetry and Perturbation Theory, Editors G. Gaeta, B. Prinari, S. Rauch-Wojciechowski, S. Terracini, World Scientific, Singapore (2005), 244--254.

\bibitem{OP1} Olver, P.J., Pohjanpelto, J., \emph{Maurer--Cartan forms and the structure of Lie pseudo-groups}, Selecta Math. {11} (2005), 99--126.

\bibitem{OP2} Olver, P.J., Pohjanpelto, J., \emph{Moving frames for Lie pseudo--groups}, Canadian J. Math. 60 (2008), 1336--1386.

\bibitem{OP3} Olver, P.J., Pohjanpelto, J., \emph{Differential invariant algebras of Lie pseudo-groups}, to appear in Adv. Math., DOI 10.1016/j.aim.2009.06.016.

\bibitem{EDS2008} Pohjanpelto, J., \emph{Reduction of exterior differential systems with infinite dimensional symmetry groups}, BIT Numerical Mathematics 48 (2008), 337--355.  

\bibitem{ShemMans2008} Shemyakova, E., Mansfield, E.L., \emph{Moving frames for Laplace invariants},
Proceedings ISSAC2008, D. Jeffrey, ed., ACM, New York, 2008, pp. 295–-302.

\bibitem{SingSter} Singer, I.M., Sternberg, S., \emph{The infinite groups of Lie and Cartan, part I, (The transitive groups)}, J. Anal. Math. 15 (1965), 1--114. 

\bibitem{Vara}
Varadarajan, V.S., \emph{Lie Groups, Lie Algebras, and Their Representations}, Springer, New York, 1984.


\end{thebibliography}
\end{document}